\documentclass[11pt]{article}
\usepackage{mathrsfs}
\usepackage{amscd}
\usepackage{amsmath,amsfonts,amssymb,amscd}
\usepackage{indentfirst,graphics,epsfig,psfrag}
\input{epsf}
\usepackage{ifpdf}
\usepackage{enumerate}
\usepackage{appendix}
\usepackage{enumerate}

\usepackage{lineno}

\makeatletter \@addtoreset{figure}{section} \makeatother
\makeatletter
\long\def\@makecaption#1#2{%
   \vskip 10\p@
   \setbox\@tempboxa\hbox{{#1}\ \ #2}%
   \ifdim \wd\@tempboxa >\hsize

       {#1}\ \ #2\par
   \else
       \hbox to\hsize{\hfil\box\@tempboxa\hfil}%
   \fi}
\makeatother

\newtheorem{thm}{Theorem}[section]
\newtheorem{cor}[thm]{Corollary}
\newtheorem{lem}[thm]{Lemma}
\newtheorem{rem}[thm]{Remark}

\newtheorem{obs}[thm]{Observation}
\newtheorem{pro}[thm]{Proposition}

\newcommand{\qed}{{\hfill\rule{3pt}{7pt}}}

\def\qed{\hfill \rule{4pt}{7pt}}

\begin{document}
\title{\bf\huge{The generalized 3-(edge) connectivity of total graphs} \footnote{Supported by NSFC No. 11561056 and 11661066.}}
\author{
\small Yinkui Li$^{1,2}$\\
\small $^{1}$Center for Applied Mathematics Tianjin University \\
\small Tianjin 300071, P.R. China\\
{\small $^{2}$College of mathematics and statistics, Qinghai Nationalities University,}\\
{\small Xining, Qinghai 810000, P.R.~China}\\[2mm]
\small E-mails:lyk463@163.com
 }
\date{}
\maketitle
\begin{abstract}
The generalized $k$-connectivity $\kappa_k(G)$ of a graph $G$,
introduced by Hager in 1985, is a natural generalization
of the concept of connectivity $\kappa(G)$, which is just for $k=2$. Total graph is generalized line graph and a large graph which obtained by incidence relation between vertices and edges of original graph. T. Hamada and T. Nonaka et al., in
\cite{Hamada} determined the connectivity of the total graph $T(G)$
for a graph $G$. In this paper, we determine the
generalized $3$-(edge)-connectivity of some total graphs and give the bounds on the
generalized $3$-(edge)-connectivity for total graph.\\[2mm]
\end{abstract}

\section{Introduction}

All graphs considered in this paper are undirected, finite and
simple. We refer to the book \cite{bondy} for graph theoretical
notation and terminology not described here. For a graph $G$, we by
$V(G)$, $E(G)$, $L(G)$, $T(G)$ denote the set of vertices, the set of edges,
the line graph and the total graph of $G$, respectively. The following we state the motivations and our
results of this paper.

Connectivity and edge-connectivity are two of the most basic
concepts of graph-theoretic subjects, both in a combinatorial sense
and an algorithmic sense. As we know, Menger's theorem is the most important basic result and fundamental theorem of connectivity. Based on this theorem, the `path' version
of \emph{connectivity} of a graph $G$ is defined as $\kappa(G)=\min\{\kappa_{G}(x,y)\,|\,x,y\in V(G),x\neq y\}$, where $\kappa_{G}(x,y)$ is the maximum number of internally disjoint paths connecting two distinct vertices $x$ and $y$ in $G$. Similarly, the \emph{edge-connectivity} of graph $G$ is defined as $\lambda(G)=\min\{\lambda_{G}(x,y)\,|\,x,y\in V(G),x\neq y\}$, where $\lambda_{G}(x,y)$ is the maximum number of edge-disjoint paths
connecting $x$ and $y$.

Although there are many elegant and powerful results on these two parameters
in graph theory, they also have their defects on measuring connection of a graph. So people want some generalizations of both
connectivity and edge-connectivity.

The generalized connectivity of a graph $G$, introduced by Hager
\cite{Hager}, is a natural and nice generalization of the `path'
version definition of connectivity. For a graph $G=(V,E)$ and a set
$S\subseteq V$ of at least two vertices, \emph{an $S$-Steiner tree}
or \emph{a Steiner tree connecting $S$} (or simply, \emph{an
$S$-tree}) is a subgraph $T=(V',E')$ of $G$ that is a tree with
$S\subseteq V'$. Two Steiner trees $T$ and $T'$ connecting $S$ are
said to be \emph{internally disjoint} if $E(T)\cap
E(T')=\varnothing$ and $V(T)\cap V(T')=S$. For $S\subseteq V(G)$ and
$|S|\geq 2$, the \emph{generalized local connectivity} $\kappa(S)$
is the maximum number of internally disjoint Steiner trees
connecting $S$ in $G$. Note that when $|S|=2$ a minimal Steiner tree
connecting $S$ is just a path connecting the two vertices of $S$.
For an integer $k$ with $2\leq k\leq n$, \emph{generalized
$k$-connectivity} (or \emph{$k$-tree-connectivity}) is defined as
$\kappa_k(G)=\min\{\kappa(S)\,|\,S\subseteq V(G),|S|=k\}$. Clearly,
when $|S|=2$, $\kappa_2(G)$ is the connectivity
$\kappa(G)$ of $G$, that is, $\kappa_2(G)=\kappa(G)$.
As a natural counterpart of the generalized connectivity, in \cite{LMS} X.Li and Y.Mao
introduced the concept of generalized
edge-connectivity, which is a generalization of the `path' version
definition of edge-connectivity. For $S\subseteq V(G)$ and $|S|\geq
2$, the \emph{generalized local edge-connectivity} $\lambda(S)$ is
the maximum number of edge-disjoint Steiner trees connecting $S$ in
$G$. For an integer $k$ with $2\leq k\leq n$, the \emph{generalized
$k$-edge-connectivity} $\lambda_k(G)$ of $G$ is then defined as
$\lambda_k(G)= \min\{\lambda(S)\,|\,S\subseteq V(G) \ and \
|S|=k\}$. It is also clear that $\lambda_2(G)=\lambda(G)$. Results on the
generalized edge-connectivity can be found in \cite{LM, LMS}.

Recently, book\cite{LM}, written by X. Li and Y. Mao, has been published, where authors bring together the known results, conjectures, and open problems on generalized connectivity and edge generalized connectivity. And thus many researchers pay more attention to the study for this topic. See \cite{Chartrand2,
LLSun, LLL2, LL, LLZ, Okamoto}.

\emph{The total graph} $T(G)$ of $G$ is graph with the vertex set
is $V(G)\cup E(G)$ and two vertices $x,y$ of $T(G)$ are adjacent if
one of the following cases hold: (i) $x,y\in V(G)$ and $x$ is
adjacent with $y$ in $G$. (ii) $x,y\in E(G)$ and $x,y$ are adjacent in
$G$. (iii) $x\in V(G)$, $y\in E(G)$, and $x,y$ are incident in $G$.
Clearly, total graph is generalized line graph.
T. Nonaka et al., in \cite{Hamada} determined the connectivity
and edge connectivity of the total graph
$T(G)$. Motivated by this, in this paper, we investigate the
generalized 3-(edge)-connectivity of total graph for some graphs, such as tree, unicycle graph, complete graph and complete bipartite graph in section 3. And further we discuss the bounds on 3-generalized-(edge-)connectivity for the total graph in section 4.

\section{Preliminary and known results}

\begin{obs}\label{obs2-1}
If $G$ is a connected graph, then $\kappa_k(G)\leq \lambda_k(G)\leq
\delta(G)$.
\end{obs}

\begin{obs}\label{obs2-2}
If $H$ is a spanning subgraph of $G$, then $\kappa_k(H)\leq
\kappa_k(G)$ and $\lambda_k(H)\leq \lambda_k(G)$.
\end{obs}

\begin{pro}\cite{Chartrand1}\label{pro2-1}
For every two integers $n$ and $k$ with $2\leq k\leq n$,
$\kappa_k(K_n)=n-\lceil k/2\rceil.$
\end{pro}

\begin{pro} \cite{LMS}\label{pro2-2}
For every two integers $n$ and $k$ with $2\leq k\leq n$,
$\lambda_k(K_n)=n-\lceil k/2\rceil.$
\end{pro}

\begin{pro} \cite{LLZ}\label{pro2-4}
Let $G$ be a connected graph of order $n\geq 6$. Then for $3\leq
k\leq 6$, $\kappa_{k}(G)\leq\kappa(G)$. Moreover, the upper bound is
always sharp for $3\leq k\leq 6$.
\end{pro}

\begin{pro} \cite{LMS}\label{pro2-5}
 For any graph $G$ of order $n$. $\lambda_{k}(G)\leq\lambda(G)$. Moreover, the upper bound is tight.
\end{pro}

\begin{pro} \cite{LLZ}\label{pro2-7}
Let $G$ be a connected graph of order $n$ with minimum degree
$\delta$. If there are two adjacent vertices of degree $\delta$,
then $\kappa_3(G)\leq \delta-1$. Moreover, the upper bound is sharp.
\end{pro}

\begin{pro} \cite{LMW}\label{pro2-8}
Let $G$ be a connected graph of order $n$ with minimum degree
$\delta$. If there are two adjacent vertices of degree $\delta$,
then $\lambda_k(G)\leq \delta-1$ for $3\leq k\leq n$. Moreover, the
upper bound is sharp.
\end{pro}

\begin{pro} \cite{LLZ}\label{pro2-9}
Let $G$ be a connected graph with $n$ vertices. For every two
integers $k$ and $r$ with $k\geq0$ and $r\in\{0,1,2,3\}$, if
$\kappa(G)=4k+r$, then $\kappa_3(G)\geq 3k+\lceil\frac{r}{2}\rceil$.
Moreover, the lower bound is sharp.
\end{pro}

\begin{thm}\cite{LLL1}\label{th2-10}
Given any three positive integers $a,b,k$ such that $a\leq b$ and
$2\leq k\leq a+b$, let $K_{a,b}$ denote a complete bipartite graph
with a bipartition of sizes $a$ and $b$. Then $\kappa_k(K_{a,b})$ is $a$ for $k\leq b-a+2$ and

$(i)$ $\frac{a+b-k+1}{2}+\Big\lfloor\frac{(a-b+k-1)(b-a+k-1)}
{4(k-1)}\Big\rfloor$ for $k>b-a+2$ and $a-b+k$ is odd.

$(ii)$ $\frac{a+b-k}{2}+\Big\lfloor\frac{(a-b+k)(b-a+k)}
{4(k-1)}\Big\rfloor$ for $k>b-a+2$ and $a-b+k$ is even.
\end{thm}

\begin{cor}\label{cor2-11}
Let $a,b$ be two integers with $2\leq a \leq b$, and $K_{a,b}$
denote a complete bipartite graph with a bipartition of sizes $a$
and $b$, respectively. Then

$$\kappa_3(K_{a,b})=\left\{
\begin{array}{ll}
a-1, &if~b=a;\\
a,&if~b>a.
\end{array}
\right.$$
\end{cor}

\begin{pro}\cite{LYMY}\label{pro2-12}
Let $a,b$ be two integers with $2\leq a \leq b$, and $K_{a,b}$
denote a complete bipartite graph with a bipartition of sizes $a$
and $b$, respectively. Then

$$\lambda_3(K_{a,b})=\left\{
\begin{array}{ll}
a-1, &if~b=a;\\
a,&if~b>a.
\end{array}
\right.$$
\end{pro}

\section{Generalized 3-(edge)-connectivity for total graph}

In this section, we determine the generalized 3-connectivity and generalized 3-edge-connectivity of some total graph $T(G)$ such as total graphs of tree, unicycle graph, complete graph and complete bipartite graph.

For $S\subset V(G)$, we by $G[S]$ denote the subgraph of $G$ which induced by $S$. Let $T$ be a subtree of $G$, then $E(T)\subset V(L(G))$. A spanning tree of induced subgraph $L(G)[E(T)]$ be called \emph{a corresponding tree} of $T$ in $L(G)$, denoted by $CT$. Now we start our investigation with tree $T_{n}$ with order $n$.

\begin{thm}\label{lem3-1}
Let $T(T_{n})$ be a total graph of tree $T_{n}$ with order $n(\geq2)$.
Then the generalized 3-connectivity of $T(T_{n})$ is
$$\kappa_{3}(T(T_{n}))=
\begin{cases}
1, &\mbox{ if $n=2$,} \\
2, &\mbox{ if $n\geq 3$.} \\
\end{cases}
$$
\end{thm}
\begin{pf}Since $T(T_{2})=T(K_{2})=C_{3}$, so $\kappa_{3}(T(T_{n}))=1$ while $n=2$.
Here we consider $n\geq 3$. Since the minimum degree of $T(T_{n})$ is
2, by Observation 2.1, we have $\kappa_{3}(T(T_{n}))\leq 2$. Next we prove
$\kappa_{3}(T(T_{n}))\geq 2$.

Suppose
$V(T_{n})=\{u_{1},u_{2},\cdots,u_{n}\}$,
$V(L(T_{n}))=\{e_{ij}|e_{ij}=u_{i}u_{j}\in E(T_{n})\}$, then $V(T(T_{n}))=V(L(T_{n}))\cup
V(T_{n})$. Let $S$ be a 3-subset of $V(T(T_{n}))$, we only need to show that there exist at least two internally disjoint $S$-trees in $T(T_{n})$.

If $|S\cap V(T_{n})|=3$, assume $S=\{u_{i},u_{j},u_{k}\}\subseteq
V(T_{n})$. By $P_{ij}$ denote the path connecting vertices $u_{i}$ and $u_{j}$ in $T_{n}$, then we obtain one $S$-tree $T=P_{ij}\cup P_{ik}$ in $T_{n}$. Let $CT$ be the corresponding tree of $T$ in $L(T_{n})$ and suppose $u_{a}, u_{b}$ and $u_{c}$ are neighbor vertices of $u_{i},u_{j}$ and $u_{k}$ in $T$, respectively. Follow this we get another $S$-tree $T'=CT\cup u_{i}e_{ia}\cup u_{j}e_{jb}\cup u_{k}e_{kc}$. Clearly, these two $S$-trees are internally disjoint, as desire. See Fig. 1 (a).

If $|S\cap V(T_{n})|=2$, assume $S=\{u_{i},u_{j},e_{ks}\}$. Then let $S'=\{u_{i},u_{j},u_{k}\}\subseteq
V(T_{n})$ and $S''=\{e_{ip},e_{jq},e_{ks}\}\subseteq
V(L(T_{n}))$. Since $T_{n}$ and $L(T_{n})$ are both connected subgraphs of in $T(T_{n})$, there exist a $S'$-tree in $T_{n}$, written as $T_{S'}$, and a $S''$-tree in $L(T_{n})$, written as $T_{S''}$. Follow this we can obtain two
internally disjoint $S$-trees:
$u_{k}e_{ks}\cup T_{S'}$ and $u_{i}e_{ip}\cup u_{j}e_{jq}\cup T_{S''}$, as desire. See Fig. 1 (b).

If $|S\cap V(T_{n})|=1$, assume $S=\{u_{i},e_{jq},e_{ks}\}$. Then let $S'=\{u_{i},u_{j},u_{k}\}\subseteq
V(T_{n})$ and  $S''=\{e_{ip},e_{jq},e_{ks}\}\subseteq
V(L(T_{n}))$. Similarly we obtain a $S'$-tree $T_{S'}$ in $T_{n}$ and a $S''$-tree $T_{S''}$ in $L(T_{n})$. Thus we construct two
internally disjoint $S$-trees in $T(T_{n})$ as: $u_{j}e_{jq}\cup u_{k}e_{ks}\cup T_{S'}$ and $u_{i}e_{ip}\cup T_{S''}$, as desire. See Fig. 1 (c).

If $|S\cap V(T_{n})|=0$, assume $S=\{e_{ip},e_{jq},e_{ks}\}$ and let $S'=\{u_{i},u_{j},u_{k}\}\subseteq
V(T_{n})$, then we can obtain a $S'$-tree in $T_{n}$, written as $T_{S'}$. Follow this we can obtain one $S$-tree in $T_{n}$ is
$u_{i}e_{ip}\cup u_{j}e_{jq}\cup u_{k}e_{ks}\cup T_{S'}$. Since $S\subset V(L(T_{n}))$ and $L(T_{n})$ be connected, we get another $S$-tree $T$ in $L(T_{n})$. Clearly, these two $S$-trees are internally disjoint, as desire. See Fig. 1 (d).

Therefore, we get $\kappa_{3}(T(T_{n}))=2$ for $n\geq 3$.\qed
\end{pf}

\begin{figure}[!hbpt]
\begin{center}
\includegraphics[scale=0.95]{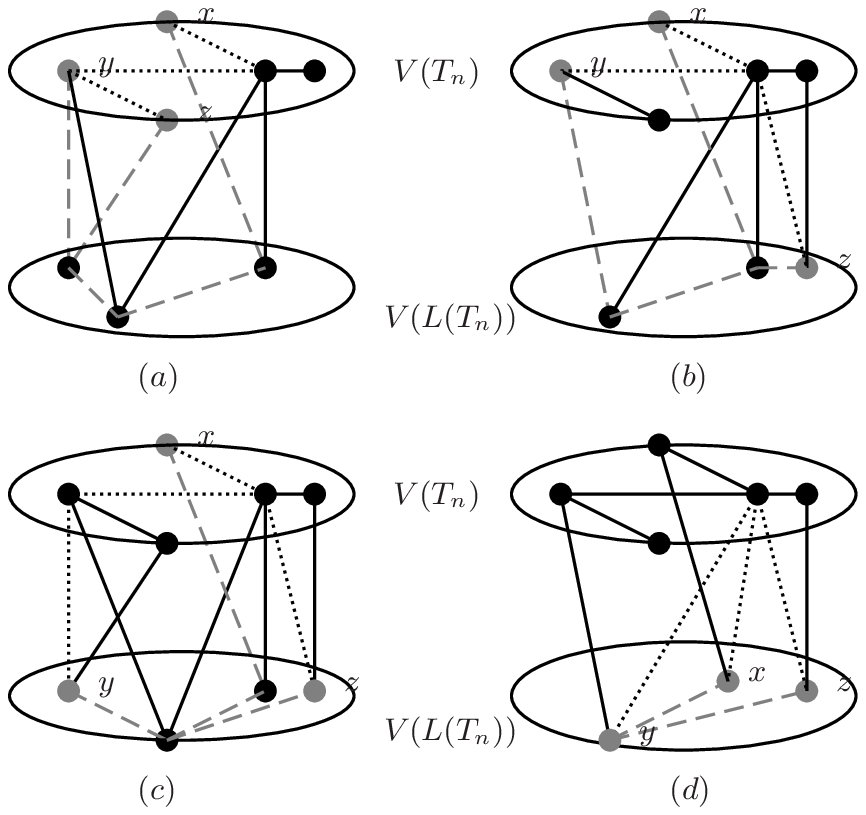}
\end{center}
\begin{center}
Figure 1: Two internally disjoint $S$-trees for 3-subset $S=\{x,y,z\}\subset V(L(T_{n}))\cup
V(T_{n})$: One is dotted, another is dashed.
\end{center}\label{fig1}
\end{figure}

\begin{thm}\label{lem3-1}
Let $T(T_{n})$ be a total graph of tree $T_{n}$ with order
$n(\geq2)$. Then the generalized 3-edge-connectivity of $T(T_{n})$
is
$$\lambda_{3}(T(T_{n}))=
\begin{cases}
1, &\mbox{ if $n=2$,} \\
2, &\mbox{ if $n\geq 3$.} \\
\end{cases}
$$
\end{thm}
\begin{pf}Since $T(T_{2})=C_{3}$, so $\lambda_{3}(T(T_{2}))=1$. While $n\geq 3$, note that the
minimum degree of $T(T_{n})$ is 2, by Observation 2.1, we get $\lambda_{3}(T(T_{n}))\leq 2$. On the other hand, by the Observation 2.1 and Theorem 3.1, $\lambda_{3}(T(T_{n}))\geq\kappa_{3}(T(T_{n}))=2$. Thus $\lambda_{3}(T(T_{n}))=2$ for $n\geq 3$.\qed
\end{pf}

The following we determine the generalized 3-connectivity and the generalized 3-edge-connectivity of unicycle graph.

\begin{thm}\label{lem3-1}
Let $T(G_{n})$ be a total graph of unicycle graph $G_{n}$ with order $n(\geq3)$.
Then the generalized 3-connectivity of $T(G_{n})$ is
$$\kappa_{3}(T(G_{n}))=
\begin{cases}
3, &\mbox{ if $G_{n}=C_{n}$,} \\
2, &\mbox{ otherwisw.} \\
\end{cases}
$$
\end{thm}

\begin{pf}
First consider $G_{n}=C_{n}$. Since $T(C_{n})$ is 4-regular graph,
by Proposition 2.7 we get $\kappa_{3}(T(C_{n}))\leq 3$. On the other
hand, since $T(C_{n})$ is 4-connected, by Proposition 2.9 we get
$\kappa_{3}(T(C_{n}))\geq 3$. Thus $\kappa_{3}(T(C_{n}))=3$.

For the general case, if $G_{n}\neq C_{n}$, since the
minimum degree of $T(G_{n})$ is 2 and thus by Observation 2.1 get
$\kappa_{3}(T(G_{n}))\leq 2$. On the other hand, suppose $C$ is unique cycle of $G_{n}$ and $e_{ij}$ is an edge of $C$. Now let $H_{n}=G_{n}-e_{ij}$, clearly, $H_{n}$ is a spanning tree of $G_{n}$ and by Theorem 3.1 we have $\kappa_{3}(T(H_{n}))=2$. Combine this with Observation 2.2 we get $\kappa_{3}(T(T_{n}))\geq 2$. Hence $\kappa_{3}(T(T_{n}))=2$ for $G_{n}\neq C_{n}$.
\qed
\end{pf}

\begin{thm}\label{lem3-1}
Let $T(G_{n})$ be a total graph of unicycle graph $G_{n}$ with order
$n(\geq3)$. Then the generalized 3-edge-connectivity of $T(G_{n})$
is
$$\lambda_{3}(T(G_{n}))=
\begin{cases}
3, &\mbox{ if $G_{n}=C_{n}$,} \\
2, &\mbox{ otherwisw.} \\
\end{cases}
$$
\end{thm}
\begin{pf} While $G_{n}=C_{n}$. Since $T(C_{n})$ is 4-regular graph, by Proposition 2.8
we get $\lambda_{3}(T(C_{n}))\leq 3$. At the same time, by Observation 2.1
and Theorem 3.3 we have $\lambda_{3}(T(C_{n}))\geq\kappa_{3}(T(C_{n}))=3$. Thus $\lambda_{3}(T(C_{n}))=3$.

While $G_{n}\neq C_{n}$, since the minimum degree of $T(G_{n})$ is 2, by Observation 2.1 we have
$\lambda_{3}(T(G_{n}))\leq 2$. On the other hand, by Observation 2.1
and Theorem 3.3 we have $\lambda_{3}(T(G_{n}))\geq\kappa_{3}(T(G_{n}))=2$, so $\lambda_{3}(T(G_{n}))=2$ for $G_{n}\neq C_{n}$. \qed
\end{pf}

Next we determine the generalized 3-connectivity and the generalized 3-edge-connectivity of complete graph. Before investigation, we first list a useful Lemma.

\begin{lem}\cite{LYMY}\label{lem3-1}
Let $L(K_{n})$ be a line graph of complete graph $K_{n}$ with $V(L(K_{n}))=\{e_{ij}|e_{ij}=u_{i}u_{j}\}$ for $V(K_{n})=\{u_{i}|1\leq i\leq n\}$. Suppose $S_{0}=\{e_{pq},e_{rs},e_{tk}\}\subseteq V(L(K_{n}))$ and $V_{S_{0}}=\{u_{p},u_{q},u_{r},u_{s},u_{t},u_{k}\}\subseteq V(K_{n})$. If the induced subgraph
$K_{n}[V_{S_{0}}]=K_{3}$, then generalized local connectivity $\kappa(S_{0})=\lfloor\frac{3(n-2)}{2}\rfloor$.
\end{lem}

\begin{lem}\cite{LYMY}
Let $K_{n}$ be complete graph with order $n(\geq3)$. Then the
generalized 3-connectivity of line graph $L(K_{n})$ is
$\kappa_{3}(L(K_{n}))=\lfloor\frac{3(n-2)}{2}\rfloor$.
\end{lem}

\begin{lem}\cite{LYMY}
Let $K_{n}$ be complete graph with order $n(\geq3)$. Then the
generalized 3-edge-connectivity of line graph $L(K_{n})$ is
$\lambda_{3}(L(K_{n}))=2n-5$.
\end{lem}

\begin{thm}\label{lem3-1}
Let $T(K_{n})$ be a total graph of complete graph $K_{n}$ with order
$n(\geq2)$. Then the generalized 3-connectivity of $T(K_{n})$ is
$$\kappa_{3}(T(K_{n}))=
\begin{cases}
3, &\mbox{ if $n=3$,} \\
\lfloor\frac{3(n-2)}{2}\rfloor+1, &\mbox{ otherwise.} \\
\end{cases}
$$
\end{thm}

\begin{pf}By Theorem 3.1 and 3.3, the result holds for cases when $n=2,3$. Now we consider $n\geq 4$. Suppose
$V(K_{n})=\{u_{1},u_{2},\cdots,u_{n}\}$ and
$V(L(K_{n}))=\{e_{ij}|e_{ij}=u_{i}u_{j}\in E(K_{n})\}$, then $V(T(K_{n}))=V(L(K_{n}))\cup
V(K_{n})$. First let $S_{0}=\{e_{ij},e_{jk},e_{ik}\}\subseteq
V(L(K_{n}))$. Clearly, $V_{S_{0}}=\{u_{i},u_{j},u_{k}\}$ and $K_{n}[V_{S_{0}}]=K_{3}$, by Lemma 3.5 there exist at most $\lfloor\frac{3(n-2)}{2}\rfloor$ internally disjoint $S_{0}$-trees in $L(K_{n})$. Besides these $S_{0}$-trees, add tree $e_{ij}u_{i}u_{j}e_{jk}\cup e_{ik}u_{i}$ together, we obtain at most $\lfloor\frac{3(n-2)}{2}\rfloor+1$ internally disjoint $S_{0}$-trees in $T(K_{n})$. Thus we get $\kappa_{3}(T(K_{n}))\leq\lfloor\frac{3(n-2)}{2}\rfloor+1$. The following we distinguish four cases to show $\kappa_{3}(T(K_{n}))\geq \lfloor\frac{3(n-2)}{2}\rfloor+1$. Let $S=\{x,y,z\}$ be a 3-subset of $V(T(K_{n}))$, we only need to show that there exist at least $\lfloor\frac{3(n-2)}{2}\rfloor+1$ internally disjoint $S$-trees in $T(K_{n})$.

\textbf{Case 1}. $|S\cap V(K_{n})|=3$

This means $x,y,z\in V(K_{n})$, assume $x=u_{a},y=u_{b},z=u_{c}$ with
$1\leq a,b,c\leq n$. Firstly, path
$zxy$ together with trees $T_{i}=u_{i}z\cup
u_{i}x\cup u_{i}y$ for $i\in \{1,2,\cdots,n\}\setminus
\{a,b,c\}$ are $n-2$ internally disjoint $S$-trees. Secondly, paths $xe_{ab}yz$, $xe_{ac}ze_{bc}y$
and trees $T_{j}=xe_{ja}e_{jb}y\cup e_{jb}e_{jc}z$ for
$j\in \{1,2,\cdots,n\}\setminus \{a,b,c\}$ are $n-1$
internally disjoint $S$-trees. Total up all we get $2n-3>\lfloor\frac{3(n-2)}{2}\rfloor+1$
internally disjoint $S$-trees in $T(K_{n})$, as desire.

\textbf{Case 2}. $|S\cap V(K_{n})|=2$

Without loss of generality, assume $x,y\in V(K_{n})$, $z\in V(L(K_{n}))$ and then let $x=u_{a},y=u_{b},z=e_{cd}$ with $1\leq
a,b,c,d\leq n$. Here first consider the case for $u_{a}u_{b}\neq u_{c}u_{d}$. If edges $u_{a}u_{b}$ and
$u_{c}u_{d}$ are nonadjacent in $K_{n}$, then we form internally disjoint
$S$-trees as: For every $i\in \{1,2,\cdots,n\}\setminus
\{a,b\}$ to form $T_{i}=xu_{i}y\cup u_{i}e_{ic}z$ and thus get $n-2$
internally disjoint $S$-trees; For every $i\in
\{1,2,\cdots,n\}\setminus \{a,b,d\}$ to form
$T'_{i}=xe_{ai}e_{bi}y\cup e_{bi}e_{id}z$ and thus get $n-3$
internally disjoint $S$-trees. Put all $T_{i}$, $T'_{i}$ with trees $yxe_{ad}z$ and
$xe_{ab}ye_{bd}z$ together, we get $2n-3$
internally disjoint $S$-trees in $T(K_{n})$.
If edges $u_{a}u_{b}$ and $u_{c}u_{d}$ are adjacent in $K_{n}$, by similar procedure as the above we also get $2n-3$
internally disjoint $S$-trees in $T(K_{n})$. Note that $2n-3>\lfloor\frac{3(n-2)}{2}\rfloor+1$, so the result holds.

Now we consider the case $u_{a}u_{b}=u_{c}u_{d}$, it is clear that for every two integers $i,j\in \{1,2,\cdots,n\}\setminus\{a,b\}$, we can get three internally disjoint $S$-trees such as $xu_{i}e_{ai}z\cup u_{i}y$, $xu_{j}e_{bj}z\cup u_{j}y$ and $xe_{aj}ze_{bi}y$ and thus get at least $\lfloor\frac{3(n-2)}{2}\rfloor$ internally disjoint $S$-trees. Put these trees with $xyz$ together, we obtain at least $\lfloor\frac{3(n-2)}{2}\rfloor+1$ internally disjoint $S$-trees in $T(K_{n})$, as we desire.

\textbf{Case 3}. $|S\cap V(K_{n})|=1$

Assume $x\in V(K_{n})$, $y,z\in V(L(K_{n}))$ and then let
$x=u_{a},y=e_{bc},z=e_{df}$ with $1\leq a,b,c,d,f\leq n$. If $e_{bc}$, $e_{df}$ and $u_{a}$ are nonadjacent each other in $T(K_{n})$, we can form internally disjoint $S$-trees as: For every $i\in \{1,2,\cdots,n\}\setminus \{a,d\}$ to form
$T_{i}=xu_{i}e_{ic}y\cup u_{i}e_{id}z$ and thus get $n-2$ internally
disjoint $S$-trees; For every $i\in
\{1,2,\cdots,n\}\setminus \{a,b\}$ to form $T'_{i}=xe_{ai}\cup
ze_{if}e_{ai}e_{ib}y$ and thus get $n-2$ internally disjoint
$S$-trees. Put all these trees with tree $ye_{ab}xu_{d}z$ together we get $2n-3>\lfloor\frac{3(n-2)}{2}\rfloor+1$ internally disjoint $S$-trees in $T(K_{n})$, as desire.

If $e_{bc}$ and $e_{df}$ are adjacent but nonadjacent to $u_{a}$, then assume $d=c$. Thus we form internally disjoint $S$-trees as: For every $i\in \{1,2,\cdots,n\}\setminus \{a,f\}$ to form $T_{i}=xu_{i}e_{ic}y\cup e_{ic}z$ and thus get $n-2$ internally
disjoint $S$-trees; For every $i\in
\{1,2,\cdots,n\}\setminus \{a,b,f\}$ to form $T'_{i}=xe_{ai}\cup
ze_{if}e_{ai}e_{ib}y$ and thus get $n-3$ internally disjoint
$S$-trees. Put all these trees with trees $xu_{f}e_{bf}y\cup e_{bf}z$ and $xe_{ab}e_{af}z\cup e_{ab}y$ together, we get
$2n-3> \lfloor\frac{3(n-2)}{2}\rfloor+1$ internally disjoint $S$-trees in $T(K_{n})$, as desire.

If $e_{bc}$ and $e_{df}$ are adjacent and one of them adjacent to $u_{a}$, then assume $d=c$ and $a=b$. Thus trees
$T_{i}=xu_{i}e_{ic}y\cup e_{ic}z$ for every $i\in \{1,2,\cdots,n\}\setminus \{f\}$ are $n-1$ internally
disjoint $S$-trees. Besides these, trees $T'_{i}=xe_{ai}\cup
ze_{if}e_{ai}e_{ib}y$ for every $i\in
\{1,2,\cdots,n\}\setminus \{a,c\}$ are also $n-2$ internally disjoint
$S$-trees. Altogether, we get
$2n-3> \lfloor\frac{3(n-2)}{2}\rfloor+1$ internally disjoint $S$-trees in $T(K_{n})$, as desire.

If $e_{bc}$, $u_{a}$ and $e_{df}$ are adjacent each other, then assume $a=d=c$. Trees $T_{i}=xu_{i}e_{ib}y\cup u_{i}e_{if}z$ for every $i\in \{1,2,\cdots,n\}\setminus \{a,f\}$ and $T'_{i}=xe_{ai}y\cup
e_{ai}z$ for every $i\in\{1,2,\cdots,n\}\setminus \{c,f\}$ are $2n-4$ internally disjoint $S$-trees. In addition to these, add tree $yu_{b}xz$ together, we get $2n-3>\lfloor\frac{3(n-2)}{2}\rfloor+1$ internally disjoint $S$-trees in $T(K_{n})$, as desire.

\textbf{Case 4}. $|S\cap V(K_{n})|=0$

It is clear $S\subseteq V(L(K_{n}))$ in this case, so by Lemma 3.6 there exist at least $\lfloor\frac{3(n-2)}{2}\rfloor$ internally disjoint $S$-trees in $V(L(K_{n}))$.
Put these $S$-trees with tree $e_{ab}u_{a}u_{c}e_{cd}\cup
u_{c}u_{g}e_{gf}$ together we get $\lfloor\frac{3(n-2)}{2}\rfloor$+1 internally disjoint $S$-trees in $T(K_{n})$, as desire.

This complete the proof.\qed
\end{pf}

\begin{thm}\label{lem3-1}
Let $T(K_{n})$ be a total graph of complete graph $K_{n}$ with order
$n(\geq2)$. Then the generalized 3-edge-connectivity of $T(K_{n})$
is $\lambda_{3}(T(K_{n}))=2n-3$.
\end{thm}
\begin{pf}Since $T(K_{n})$ is $2n-2$ regular graph, by proposition
2.8 $\lambda_{3}(T(K_{n}))\leq 2n-3$. Next we by constructing $2n-3$ edge disjoint $S$-trees in $T(K_{n})$ for any 3-subset $S$ of $V(T(K_{n}))$ to prove $\lambda_{3}(T(K_{n}))\geq2n-3$.

Recall of the proof of Theorem 3.8, except case
$|S\cap V(K_{n})|=0$ and case $|S\cap V(K_{n})|=2$ for $S=\{u_{a},u_{b},e_{cd}\}$ with $u_{a}u_{b}=u_{c}u_{d}$, there always exist at least $2n-3$ internally disjoint $S$-trees in $T(K_{n})$, which are also edge disjoint $S$-trees in $T(K_{n})$. Thus here we only need to consider the above two exception cases and show there still exist at least $2n-3$ edge disjoint $S$-trees in $T(K_{n})$.

If $|S\cap V(K_{n})|=0$, assume $S=\{e_{ab},e_{cd},e_{gf}\}\subset V(L(K_{n}))$, by Lemma 3.7 there exist at least $2n-5$ edge disjoint $S$-trees in $L(K_{n})$, in addition these trees, add two edge disjoint $S$-trees $e_{ab}u_{a}u_{c}e_{cd}\cup
u_{c}u_{g}e_{gf}$ and $e_{ab}u_{b}u_{d}e_{cd}\cup
u_{d}u_{f}e_{gf}$ together, we obtain at least $2n-3$ edge disjoint $S$-trees in $T(K_{n})$, as desire.

If $|S\cap V(K_{n})|=2$ for $S=\{u_{a},u_{b},e_{cd}\}$ with $u_{a}u_{b}=u_{c}u_{d}$, this means $S=\{u_{a},u_{b},e_{ab}\}$, then we form $S$-trees $T_{i}^{1}=xu_{i}e_{bi}z\cup e_{bi}y$ and $T_{i}^{2}=xe_{ai}u_{i}y\cup e_{ai}z$ for every $i\in
[n]\setminus \{ a,b\}$ and thus get $2n-4$ edge disjoint $S$-trees. Add tree $xzy$ together with all $T_{i}^{1}$ and $T_{i}^{2}$ we get $2n-3$ edge disjoint $S$-trees in $T(K_{n})$, as desire.

Thus we get $\lambda_{3}(T(K_{n}))\geq2n-3$. This complete the proof.\qed
\end{pf}

At end of this section, we determine the generalized 3-connectivity and the generalized 3-edge-connectivity of the total graph of the complete bipartite graph. We start with definition and Lemmas.

\emph{The Cartesian product} $G_{1}\times G_{2}$ of $G_{1}$ and $G_{2}$ is a graph which has vertex
set $V(G_{1})\times V(G_{2})$ with two vertices $x=(u,u')$ and $y=(v,v')$ adjacent iff for $u=v$, $u'$ is adjacent with $v'$ in $G_{2}$ or
$u'=v'$, $u$ is adjacent with $v$ in $G_{1}$.

It is clear that line graph $L(K_{m,n})$ of complete bipartite graph $K_{m,n}$ is the Cartesian product of $K_{m}$ and $K_{n}$ and the generalized 3-connectivity of $L(K_{m,n})$ has been determined in our another paper, which listed as follow.

\begin{lem}\cite{LYMY}
Let $L(K_{m,n})$ be the line graph of complete bipartite graph $K_{m,n}(m\leq n)$, then the generalized
3-connectivity of $L(K_{m,n})$ is
$\kappa_{3}(L(K_{m,n}))=\kappa_{3}(K_m\times K_n)=m+n-3$.
\end{lem}

\begin{thm}\label{lem3-1}
Let $T(K_{m,n})$ be a total graph of complete bipartite graph
$K_{m,n}$ with $1\leq m\leq n$. Then the generalized
3-connectivity of $T(K_{m,n})$ is
$$\kappa_{3}(T(K_{m,n}))=
\begin{cases}
2m-1, &\mbox{ if  $m=n$,} \\
2m,   &\mbox{ if  $m<n$.} \\
\end{cases}
$$
\end{thm}

\begin{pf} Suppose $U=\{u_1,u_2,\cdots,u_m\}$ and $V=\{v_1,v_2,\cdots,v_{n}\}$ be
the two parts of $K_{m,n}$, ie., $V(K_{m,n})=U\cup V$. Then $V(L(K_{m,n}))=\{e_{ij}|e_{ij}=u_{i}v_{j}\}$ and $V(T(K_{m,n}))=U\cup V\cup V(L(K_{m,n}))$. Clearly, by Theorem 3.1, the results hold for case $m=1$. The following we consider cases for $n\geq m\geq 2$.

\textbf{Case 1}. $m=n$

Since $T(K_{m,m})$ is $2m$ regular, by proposition 2.7 we get $\kappa_{3}(T(K_{m,m}))\leq 2m-1$. The following we prove $\kappa_{3}(T(K_{m,m}))\geq 2m-1$. In fact, for any 3-subset $S=\{x,y,z\}$ of $V(T(K_{m,m}))$, here we only need to show that there exist at least $2m-1$ internally disjoint $S$-trees in $T(K_{m,m})$.

\textbf{Subcase 1.1}. $|S\cap V(K_{m,m})|=3$

$|S\cap V(K_{m,m})|=3$ means $S\subseteq V(K_{m,m})$. Then either $|S\cap V|=3$ or $|S\cap V|=2$ and $|S\cap U|=1$. If $|S\cap V|=3$, assume $x=v_{1},y=v_{2},z=v_{3}$ (see Fig 2 (a)), we can form $2m(>2m-1)$ internally disjoint $S$-trees such as: $T_{i}=xu_{i}z\cup u_{i}y$ and $T'_{i}=xe_{i1}e_{i2}y\cup e_{i2}e_{i3}z$ for $1\leq i\leq m$. If $|S\cap V|=2$ and $|S\cap U|=1$, assume $x=v_{1},y=v_{2},z=u_{m}$ (see Fig 2 (b)), we form $2m-1$ internally disjoint $S$-trees as follow: $xzy$, $xe_{m1}e_{11}e_{12}y\cup e_{m1}z$, $xu_{1}ye_{m2}z$ and $xu_{i}y\cup u_{i}v_{i+1}z$, $xe_{i1}e_{i2}e_{i(i+1)}e_{m(i+1)}z\cup e_{i2}y$ for $2\leq i\leq m-1$. All as we desire.
\begin{figure}[!hbpt]
\begin{center}
\includegraphics[scale=0.95]{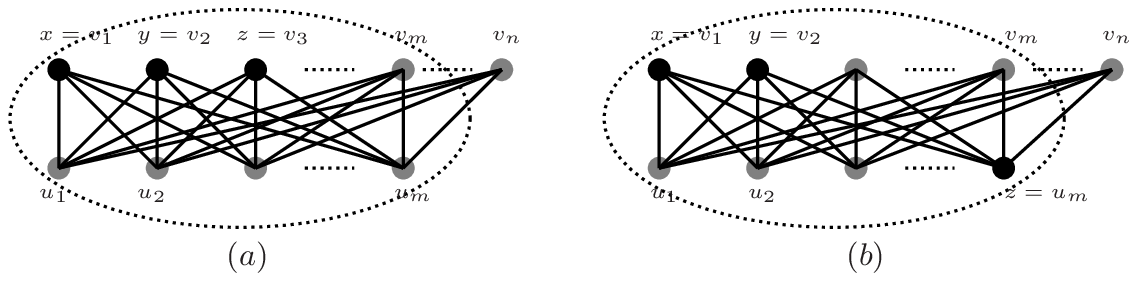}
\end{center}
\begin{center}
Figure 2: The subcase for $|S\cap V(K_{m,n})|=3$.
\end{center}\label{fig2}
\end{figure}

\textbf{Subcase 1.2}. $|S\cap V(K_{m,m})|=2$

Here we need discuss five possible cases. If $x,y,z$ are nonadjacent each other in $T(K_{m,m})$, assume $x=v_{1},y=v_{2},z=e_{mm}(m\neq 1,2)$ (see Fig 3 (a)). Then $2m-1$ internally disjoint $S$-trees be formed as follows: $xu_{m}z\cup u_{m}y$, $xu_{1}v_{m}z\cup u_{1}y$, $xe_{m1}e_{m2}z\cup e_{m2}y$, $xe_{11}e_{12}e_{1m}z\cup e_{12}y$ and $xu_{i}e_{im}z\cup u_{i}y$ for $2\leq i\leq m-1$, $xe_{i1}e_{i2}e_{i(i+1)}e_{m(i+1)}z\cup e_{i2}y$ for $2\leq i\leq m-2$.

If $x$ is nonadjacent to $y,z$ but $y$ and $z$ are adjacent in $T(K_{m,m})$, assume $x=v_{1},y=v_{2},z=e_{22}$ (see Fig 3 (b)). Then $2m-1$ internally disjoint $S$-trees are formed as: $xu_{2}z\cup u_{2}y$, $xu_{1}yz$, $xe_{11}e_{12}z\cup e_{12}y$ and $xu_{i}v_{i}e_{2i}z\cup u_{i}y$, $xe_{i1}e_{i2}z\cup e_{i2}y$ for $3\leq i\leq m$.

If $z$ is nonadjacent to $x,y$ but $x$ and $y$ are adjacent in $T(K_{m,m})$, assume $x=v_{1},y=u_{1},z=e_{mm}$ (see Fig 3 (c)). Then $2m-1$ internally disjoint $S$-trees be formed as: $xyv_{m}z$, $ye_{11}xu_{m}z$, $xe_{m1}ze_{1m}y$ and $yv_{i}u_{i}x\cup u_{i}e_{im}z$, $xe_{i1}e_{ii}e_{1i}y\cup e_{1i}e_{mi}z$ for $2\leq i\leq m-1$.

If $x$ and $z$ are both adjacent to $y$ but $x$ and $z$ are nonadjacent in $T(K_{m,m})$, assume $x=v_{1},y=u_{1},z=e_{1m}$ (see Fig 3 (d)). Then $2m-1$ internally disjoint $S$-trees be formed as: $xyz$, $xe_{11}y\cup e_{11}z$, $xu_{m}v_{m}z\cup u_{m}y$ and $xu_{i}v_{i}y\cup u_{i}e_{im}z$, $xe_{i1}e_{ii}e_{1i}y\cup e_{ii}z$ for $2\leq i\leq m-1$.

If $x,y,z$ are adjacent each other in $T(K_{m,m})$, assume $x=v_{1},y=u_{1},z=e_{11}$ (see Fig 3 (e)). Then $xzy$, $xe_{1m}z\cup e_{1m}e_{mm}u_{m}y$, $ye_{m1}z\cup e_{m1}e_{m2}v_{2}x$ and $xe_{1i}z\cup e_{1i}e_{ii}u_{i}y$, $ye_{i1}e_{i(i+1)}v_{i+1}x\cup e_{i1}z$ for $2\leq i\leq m-1$ are $2m-1$ internally disjoint $S$-trees in $T(K_{m,m})$.

\begin{figure}[!hbpt]
\begin{center}
\includegraphics[scale=0.95]{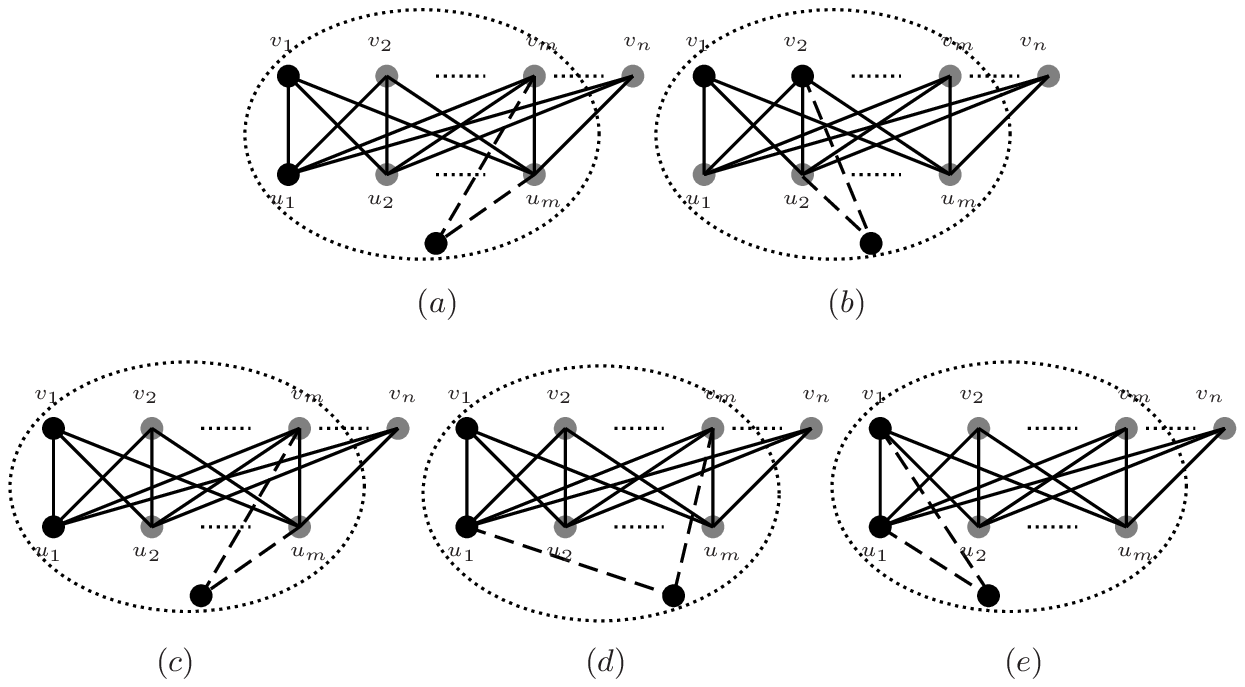}
\end{center}
\begin{center}
Figure 3: The subcases for $|S\cap V(K_{m,n})|=2$.
\end{center}\label{fig3}
\end{figure}

\textbf{Subcase 1.3}. $|S\cap V(K_{m,m})|=1$

Similarly, we also need discuss five possible cases while $|S\cap V(K_{m,m})|=1$. If $x,y,z$ are nonadjacent each other in $T(K_{m,m})$, assume $x=v_{1},y=e_{22},z=e_{mm}$ (see Fig 4 (a)). Then $2m-1$ internally disjoint $S$-trees be formed as: $xu_{2}v_{m}z\cup u_{2}y$, $xu_{m}v_{2}y\cup u_{m}z$, $xu_{1}e_{12}y\cup e_{12}e_{1m}z$, $xe_{21}ye_{2m}z$, $xe_{m1}ze_{m2}y$ and $xu_{i}v_{i}e_{mi}z\cup v_{i}e_{2i}y$, $xe_{i1}e_{ii}e_{im}z\cup e_{ii}e_{i2}y$ for $3\leq i\leq m-1$.

If $y,z$ are adjacent but they are nonadjacent to $x$ in $T(K_{m,m})$, assume $x=v_{1},y=e_{22},z=e_{m2}$ (see Fig 4 (b)). Then $2m-1$ internally disjoint $S$-trees be constructed as: $xu_{1}v_{2}z\cup v_{2}y$, $xe_{11}e_{12}z\cup e_{12}y$, $xu_{2}v_{m}u_{m}z\cup u_{2}y$, $xe_{21}y\cup e_{21}e_{2m}e_{mm}z$, $xe_{m1}zy$ and $xu_{i}e_{i2}y\cup e_{i2}z$, $xe_{i1}e_{ii}e_{2i}y\cup e_{ii}e_{mi}z$ for $3\leq i\leq m-1$.

If $x,y$ are adjacent but they are nonadjacent to $z$ in $T(K_{m,m})$, assume $x=v_{1},y=e_{21},z=e_{mm}$ (see Fig 4 (c)). Then $2m-1$ internally disjoint $S$-trees be constructed as: $yxu_{m}z$, $xu_{2}v_{m}z\cup u_{2}y$, $xe_{m1}ze_{2m}y$, $xu_{1}v_{2}e_{22}y\cup e_{22}e_{m2}z$, $xe_{11}y\cup e_{11}e_{1m}z$ and $xu_{i}v_{i}e_{2i}y\cup v_{i}e_{mi}z$, $xe_{i1}y\cup e_{i1}e_{im}z$ for $3\leq i\leq m-1$.

If $x$ and $z$ are both adjacent to $y$ but $x$ and $z$ are nonadjacent in $T(K_{m,m})$, assume $x=v_{1},y=e_{21},z=e_{22}$ (see Fig 4 (d)). Then $2m-1$ internally disjoint $S$-trees be constructed as: $xyz$, $xu_{2}y\cup u_{2}z$, $xe_{11}u_{1}v_{2}z\cup e_{11}y$ and $xu_{i}v_{i}e_{2i}y\cup e_{2i}z$, $xe_{i1}y\cup e_{i1}e_{i2}z$ for $3\leq i\leq m$.

If $x,y,z$ are adjacent each other in $T(K_{m,m})$, assume $x=v_{1},y=e_{11},z=e_{21}$ (see Fig 4 (e)). Then $2m-1$ internally disjoint $S$-trees be constructed as: $yxz$, $xu_{1}y\cup u_{1}v_{2}u_{2}z$, $zye_{12}e_{22}u_{2}x$ and $xu_{i}v_{i}e_{2i}z\cup v_{i}e_{1i}y$, $xe_{i1}y\cup e_{i1}z$ for $3\leq i\leq m$.

\begin{figure}[!hbpt]
\begin{center}
\includegraphics[scale=0.95]{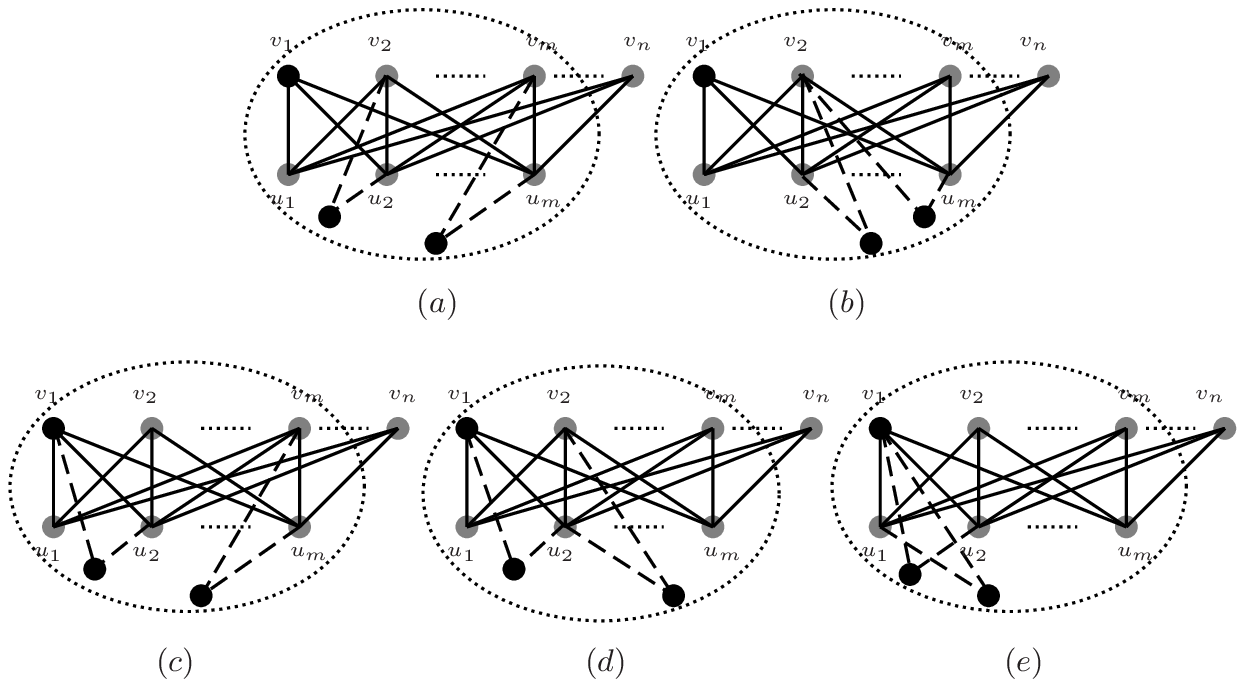}
\end{center}
\begin{center}
Figure 4: The subcases for $|S\cap V(K_{m,n})|=1$.
\end{center}\label{fig5}
\end{figure}

\textbf{Subcase 1.4}. $|S\cap V(K_{m,m})|=0$

$|S\cap V(K_{m,m})|=0$ means $S\subseteq V(L(K_{m,m}))$, without loss generality, suppose $S=\{e_{pq},e_{rs},e_{tk}\}$, then by Lemma 3.10 there always exist $2m-3$ internally disjoint $S$-trees in $L(K_{m,m})$, named as $T_{i}(1\leq i\leq 2m-3)$. Now we add two $S$-trees $e_{pq}u_{p}v_{j}u_{t}e_{tk}\cup v_{j}u_{r}e_{rs}$ and $e_{pq}v_{q}u_{j}v_{k}e_{tk}\cup u_{j}v_{s}e_{rs}$ to all $T_{i}$, here $j\neq p,q,r,s,t,k$. In total, we get $2m-1$ internally disjoint $S$-trees in $T(K_{m,m})$, as desire.

\textbf{Case 2}. $m<n$

Since the minimum degree of $T(K_{m,n})$ is $2m$, by Observation 2.1 we get
$\kappa_{3}(T(K_{m,n}))\leq 2m$. Let $S=\{x,y,z\}$ be a 3-subset of $V(T(K_{m,n}))$, now we by forming $2m$ internally disjoint $S$-trees in $T(K_{m,n})$ to prove $\kappa_{3}(T(K_{m,n}))\geq 2m$.

\textbf{Subcase 2.1}. $|S\cap V(K_{m,n})|=3$

Clearly, either $|S\cap V|=3$ or $|S\cap V|=2$ and $|S\cap U|=1$. If $|S\cap V|=3$, by the proof of Subcase 1.1, there exist $2m$ internally disjoint $S$-trees in $T(K_{m,n})$. Here we consider the case for $|S\cap V|=2$ and $|S\cap U|=1$, assume $x=v_{1},y=v_{2},z=u_{m}$ (see Fig 2 (b)),we can form $2m$ internally disjoint $S$-trees as: $xzy$, $xe_{m1}ze_{m2}y$ and $xu_{i}y\cup u_{i}v_{i+2}z$, $xe_{i1}e_{i2}e_{i(i+2)}e_{m(i+2)}z\cup e_{i2}y$ for $1\leq i\leq m-1$. As desire.

\textbf{Subcase 2.2}. $|S\cap V(K_{m,n})|=2$

Here we need discuss five possible cases. If $x,y,z$ are nonadjacent each other in $T(K_{m,n})$, assume $x=v_{1},y=v_{2},z=e_{mm}$ (see Fig 3 (a)). Then we form $2m$ internally disjoint $S$-trees as: $xu_{m}z\cup u_{m}y$, $xe_{m1}e_{m2}y\cup e_{m2}z$ and $xu_{i}v_{i+2}e_{m(i+2)}z\cup u_{i}y$, $ye_{i2}e_{i1}x\cup e_{i1}e_{im}z$ for $1\leq i\leq m-1$.

If $y$ and $z$ are nonadjacent to $x$ but $y$ and $z$ are adjacent in $T(K_{m,n})$, assume $x=v_{1},y=v_{2},z=e_{11}$ (see Fig 3  (b)). Then $2m$ internally disjoint $S$-trees be constructed as: $xu_{1}z\cup u_{1}y$, $ye_{12}zx$ and $xu_{i}v_{i+1}e_{1(i+1)}z\cup u_{i}y$, $xe_{i1}e_{i2}y\cup e_{i1}z$ for $2\leq i\leq m$.

If $x$ and $y$ are nonadjacent to $z$ but $x$ and $y$ are adjacent in $T(K_{m,n})$, assume $x=v_{1},y=u_{1},z=e_{mm}$ (see Fig 3  (c)). Then $2m$ internally disjoint $S$-trees be constructed as: $xyv_{m}z$, $yv_{m+1}u_{m}z\cup u_{m}x$, $xe_{11}y\cup e_{11}e_{m1}z$, $xe_{m1}ze_{m(m+1)}e_{1(m+1)}y$ and $yv_{i}u_{i}x\cup u_{i}e_{im}z$, $xe_{i1}e_{ii}e_{1i}y\cup e_{1i}e_{mi}z$ for $2\leq i\leq m-1$.

If $x$ and $z$ are both adjacent to $y$ but $x$ and $z$ are nonadjacent in $T(K_{m,n})$, assume $x=v_{1},y=u_{1},z=e_{1m}$ (see Fig 3 (d)). Then $2m$ internally disjoint $S$-trees be constructed as follows: $xyz$, $xe_{11}y\cup e_{11}z$, $xu_{m}v_{m}z\cup u_{m}y$, $xe_{m1}e_{m(m+1)}e_{1(m+1)}y\cup e_{1(m+1)}z$ and $xu_{i}v_{i}y\cup u_{i}e_{im}z$, $xe_{i1}e_{ii}e_{1i}y\cup e_{ii}z$ for $2\leq i\leq m-1$.

If $x,y,z$ are adjacent each other in $T(K_{m,n})$, assume $x=v_{1},y=u_{1},z=e_{11}$ (see Fig 3 (e)). Then $xzy$, $xy\cup xe_{1(m+1)}z$, $xe_{1m}z\cup e_{1m}e_{mm}u_{m}y$, $ye_{m1}z\cup e_{m1}e_{m2}v_{2}x$ and $xe_{1i}z\cup e_{1i}e_{ii}u_{i}y$, $ye_{i1}e_{i(i+1)}v_{i+1}x\cup e_{i1}z$ for $2\leq i\leq m-1$ are $2m$ internally disjoint $S$-trees in $T(K_{m,n})$, as desire.

\textbf{Subcase 2.3}. $|S\cap V(K_{m,n})|=1$

Similarly, we also need discuss five possible cases. If $x,y,z$ are nonadjacent each other in $T(K_{m,n})$, assume $x=v_{1},y=e_{22},z=e_{mm}$ (see Fig 4 (a)). Then $2m$ internally disjoint $S$-trees be constructed as: $xu_{2}v_{m}z\cup u_{2}y$, $xu_{m}v_{2}y\cup u_{m}z$, $xu_{1}e_{12}y\cup e_{12}e_{1m}z$, $xe_{21}ye_{2m}z$, $xe_{m1}ze_{m2}y$, $xe_{11}e_{1(m+1)}e_{2(m+1)}y\cup e_{2(m+1)}e_{m(m+1)}z$ and $xu_{i}v_{i}e_{mi}z\cup e_{mi}e_{2i}y$, $xe_{i1}e_{im}z\cup e_{im}e_{i2}y$ for $3\leq i\leq m-1$.

If $y$ and $z$ are adjacent but they are nonadjacent to $x$ in $T(K_{m,n})$, assume $x=v_{1},y=e_{22},z=e_{m2}$ (see Fig 4 (b)). Then $2m$ internally disjoint $S$-trees be constructed as: $xu_{1}v_{2}z\cup v_{2}y$, $xe_{11}e_{12}z\cup e_{12}y$, $xu_{2}v_{m+1}e_{m(m+1)}z\cup u_{2}y$, $xu_{m}z\cap u_{m}e_{m2}y$, $xe_{21}yz$, $xe_{m1}z\cup e_{m1}e_{mm}e_{2m}y$ and $xu_{i}e_{i2}y\cup e_{i2}z$, $xe_{i1}e_{ii}e_{2i}y\cup e_{ii}e_{mi}z$ for $3\leq i\leq m-1$.

If $x$ and $y$ are adjacent but they are nonadjacent to $z$ in $T(K_{m,n})$, assume $x=v_{1},y=e_{11},z=e_{mm}$ (see Fig 4 (c)). Then $2m$ internally disjoint $S$-trees be constructed as: $xye_{1m}z$, $xu_{1}v_{m}z\cup u_{1}y$, $xe_{m1}z\cup e_{1m}y$, $xu_{m}e_{m(m+1)}e_{1(m+1)}y\cup e_{m(m+1)}z$ and $xu_{i}v_{i}e_{1i}y\cup v_{i}e_{mi}z$, $xe_{i1}y\cup e_{i1}e_{im}z$ for $2\leq i\leq m-1$.

If $x$ and $z$ are both adjacent to $y$ but $x$ and $z$ are nonadjacent in $T(K_{m,n})$, assume $x=v_{1},y=e_{21},z=e_{22}$ (see Fig 4 (d)). Then $2m$ internally disjoint $S$-trees be constructed as: $xyz$, $xu_{2}y\cup u_{2}z$, $xe_{11}e_{12}z\cup e_{11}y$, $xu_{1}v_{m+1}e_{2(m+1)}y\cup e_{2(m+1)}z$ and $xu_{i}v_{i}e_{2i}y\cup e_{2i}z$, $xe_{i1}y\cup e_{i1}e_{i2}z$ for $3\leq i\leq m$.

If $x,y,z$ are adjacent each other in $T(K_{m,n})$, assume $x=v_{1},y=e_{11},z=e_{21}$ (see Fig 4 (e)). Then $2m$ internally disjoint $S$-trees be constructed as: $yxz$, $xu_{1}y\cup u_{1}e_{12}e_{22}z$, $xe_{21}z\cup e_{21}y$, $xu_{2}e_{2(m+1)}z\cup e_{2(m+1)}e_{1(m+1)}y$ and $xu_{i}v_{i}e_{2i}z\cup v_{i}e_{1i}y$, $xe_{i1}y\cup e_{i1}z$ for $3\leq i\leq m$.

\textbf{Subcase 2.4}. $|S\cap V(K_{m,n})|=0$

$|S\cap V(K_{m,n})|=0$ means $S\subseteq V(L(K_{m,n}))$, assume $S=\{e_{pq},e_{rs},e_{tk}\}$, then by Lemma 3.10 there exist $m+n-3$ internally disjoint $S$-trees in $L(K_{m,n})$, named as $T_{i}(1\leq i\leq m+n-3)$. Now we add two $S$-trees $e_{pq}u_{p}v_{j}u_{t}e_{tk}\cup v_{j}u_{r}e_{rs}$ and $e_{pq}v_{q}u_{j}v_{k}e_{tk}\cup u_{j}v_{s}e_{rs}$ to all $T_{i}$, here $j\neq p,q,r,s,t,k$. Thus we get $m+n-1(\geq2m)$ internally disjoint $S$-trees in $T(K_{m,n})$, as desire.

By now we complete the proof.\qed
\end{pf}

Note that the fact $T(K_{m,n})$ is $2m$
-regular graph while $m=n$ and its minimum degree is $2m$
while $m<n$. Combine this with Observation 2.1 and Theorem 3.11, we immediately get Theorem 3.12.
\begin{thm}
Let $T(K_{m,n})$ be a total graph of complete bipartite graph
$K_{m,n}(m\leq n)$. Then the generalized
3-edge-connectivity of $T(K_{m,n})$ is
 $$\lambda_{3}(T(K_{m,n}))=
\begin{cases}
2m-1, &\mbox{ if $m=n$,} \\
2m, &\mbox{ if $m<n$.} \\
\end{cases}
$$
\end{thm}

\section{Bound for generalized 3-connectivity of total graph}

In this section we give some bounds for the generalized
3-connectivity of total graph $T(G)$.

\begin{lem}\cite{HAMADA}
Let $G$ be graph for which $\kappa(G)\geq m$. Then $\kappa(T(G))\geq
2m$ and $\lambda(T(G))\geq2m$.
\end{lem}

\begin{lem}\cite{LYMY}
Let $G$ be connected graph, then
$\kappa_{3}(L(G))\geq\lambda_{3}(G)$.
\end{lem}

\begin{thm}
Let $G$ be connected graph with connectivity $\kappa(G)$ and minimum degree $\delta(G)$. Then
$$\lfloor\frac{3\kappa(G)-1}{2}\rfloor\leq\kappa_{3}(T(G))\leq2\delta(G)$$
\end{thm}

\begin{pf} Since the minimum degree of $T(G)$ is $2\delta(G)$, so by Proposition 2.6 we get
$\kappa_{3}(T(G))\leq 2\delta(G)$. On the other hand,
by Lemma 4.1, we know $\kappa(T(G))\geq 2\kappa(G)$ and let $\kappa(T(G))=4a+b$ with $b\in\{0,1,2,3\}$. Then we get $a\geq\frac{2\kappa(G)-b}{4}$. Combine this with
Proposition 2.9 we have
$$\kappa_{3}(T(G))\geq 3a+\lceil \frac{b}{2}\rceil
\geq \frac{3(2\kappa(G)-b)}{4}+\lceil \frac{b}{2}\rceil$$

Since $b\in\{0,1,2,3\}$, so the value of $\frac{3(2\kappa(G)-b)}{4}+\lceil \frac{b}{2}\rceil$ can meet $\frac{3\kappa(G)}{2}, \frac{3\kappa(G)}{2}+\frac{1}{4}, \frac{3\kappa(G)}{2}-\frac{1}{4}$ and $\frac{3\kappa(G)}{2}-\frac{1}{2}$, respectively. Consider this with $\kappa_{3}(T(G))$ is an integer, we get
$$\kappa_{3}(T(G))\geq\lfloor\frac{3\kappa(G)-1}{2}\rfloor$$\qed
\end{pf}

\begin{rem}
The upper bound is sharp for complete bipartite graph $K_{m,n}$ with $n>m$. The lower bound is also sharp for complete graph $K_{n}$ with $n\neq 3$.
\end{rem}

Similar we also get the following result for the generalized
3-edge-connectivity of total graph $T(G)$.

\begin{thm}
Let $G$ be connected graph with edge connectivity $\lambda(G)$ and minimum degree $\delta(G)$. Then
$$\min\{2\lambda(G)-1, 2\lambda_{3}(G), \lambda_{3}(G)+2\}\leq\lambda_{3}(T(G))\leq2\delta(G)$$
\end{thm}

\begin{pf} Since the minimum degree of $T(G)$ is $2\delta(G)$, by Observation 2.1 we get
$\lambda_{3}(T(G))\leq 2\delta(G)$. Next we prove $\lambda_{3}(T(G))\geq \min\{2\lambda(G)-1, 2\lambda_{3}(G), \lambda_{3}(G)+2\}$.

Suppose $V(G)=\{v_1,v_2,\cdots,v_{n}\}$ and $V(L(G))=\{e_{ij}|v_{i}v_{j}\in E(G)\}$, let $S=\{x,y,z\}$ be a 3-subset of $V(T(G))$. The following we by constructing edge disjoint $S$-trees in $T(G)$ to prove $\lambda_{3}(T(G))\geq \min\{2\lambda(G)-1, 2\lambda_{3}(G), \lambda_{3}(G)+2\}$. Now distinguish four cases to complete the proof.

\textbf{Case 1}. $|S\cap V(G)|=3$

$|S\cap V(G)|=3$ means $S\subseteq V(G)$, since there are $\lambda_{3}(G)$ edge disjoint $S$-trees in $G$, name them as $T_{i}$. In addition to these, it is clear that every $T_{i}$ has a corresponding tree $CT_{i}$ in $L(G)$ and every $CT_{i}$ can be formed as a $S$-tree by connecting it with $S$. Thus altogether we can get at least $2\lambda_{3}(G)$ edge disjoint $S$-trees in $T(G)$.

\textbf{Case 2}. $|S\cap V(G)|=2$

Suppose $x,y\in V(G)$, $z\in V(L(G))$ and let $x=v_{i},y=v_{j},z=e_{pq}$. If $T(G)[S]$ is a triangle, since there are $\lambda(G)$ edge disjoint $xy$-paths in $G$ and thus assume $xv_{1}v_{2}\cdots v_{k}y$ is a $xy$-path. Based on this $xy$-path we can form 2 edge disjoint $S$-trees in $T(G)$ as: $xe_{i1}v_{1}v_{2}\cdots v_{k}y\cup e_{i1}z$ and $x v_{1}e_{12}e_{23}\cdots e_{kj}y\cup e_{kj}z$. Thus we total get $2\lambda(G)$ edge disjoint $S$-trees in $T(G)$.

If $p=i$ but $q\neq j$, assume $S'=\{v_{i}, v_{j}, v_{q}\}\subseteq V(G)$. Since there exist $\lambda_{3}(G)$ edge disjoint $S'$-trees in $G$, name them as $T_{i}$. For every pair trees $T_{i}$ and its corresponding tree $CT_{i}$, we by symmetric difference operation on $G$ and $L(G)$ can obtain $2\lambda_{3}(G)$ edge disjoint $S$-trees in $T(G)$.

If $p\neq i$ and $q\neq j$, assume $S'=\{v_{i}, v_{j}, v_{q}\}$. Since there exist $\lambda_{3}(G)$ edge disjoint $S'$-trees in $G$, named as $T_{i}$. For every pair trees $T_{i}$ and its corresponding tree $CT_{i}$, we also can obtain at least$\lambda_{3}(G)+2$ edge disjoint $S$-trees in $T(G)$.

\textbf{Case 3}. $|S\cap V(G)|=1$

Assume $x\in V(G)$, $y,z\in V(L(G))$ and let $x=v_{i},y=e_{jk},z=e_{pq}$. No matter what vertices $x,y,z$ they are, there exist 3 elements among $v_{i},v_{j},v_{k},v_{p},v_{q}$ which can form a 3-subset $S'\subseteq V(G)$. Since there exist $\lambda_{3}(G)$ edge disjoint $S'$-trees in $G$, named as $T_{i}$. Similarly, based on every pair $T_{i}$ and its corresponding tree $CT_{i}$ we get at least $\lambda_{3}(G)+2$ edge disjoint $S$-trees in $T(G)$.

\textbf{Case 4}. $|S\cap V(G)|=0$

This case means $S\subseteq V(L(G))$, by Lemma 4.2 there exist $\lambda_{3}(G)$ internally disjoint $S$-trees in $L(G)$, which are also edge disjoint $S$-trees. Based on this, using vertices in graph $G$ can obtain at least $\lambda_{3}(G)+2$ edge disjoint $S$-trees in $T(G)$.

By argument of the above , we can claim that there exist at least $\min\{2\lambda(G)-1, 2\lambda_{3}(G), \lambda_{3}(G)+2\}$ edge disjoint $S$-trees in $T(G)$. Thus $\lambda_{3}(T(G))\geq \min\{2\lambda(G)-1, 2\lambda_{3}(G), \lambda_{3}(G)+2\}$.\qed
\end{pf}

\begin{rem}
The upper bound is sharp for complete bipartite graph $K_{m,n}$ with $n>m$. The lower bound is also sharp for $\lambda_{3}(T(C_{n}))=2\lambda(C_{n})-1=\lambda_{3}(C_{n})+2$, $\lambda_{3}(T(K_{n}))=2\lambda(K_{n})-1$ and $\lambda_{3}(T(K_{m,n}))=2\lambda_{3}(K_{m,n})$ with $n>m$.
\end{rem}


\begin{thebibliography}{1}

\bibitem{bondy} J.A. Bondy, U.S.R. Murty, \emph{Graph Theory},
GTM 244, Springer, 2008.

\bibitem{Beineke} L.W. Beineke, R.J. Wilson, \emph{Topics in Structural Graph
Theory}, Cambrige University Press, 2013.

\bibitem{Boesch} F.T. Boesch, S. Chen, \emph{A generalization of line connectivity
and optimally invulnerable graphs}, SIAM J. Appl. Math. 34(1978),
657--665.

\bibitem{Chang} S. Chang, \emph{The uniqueness and nonuniqueness of the
triangular association scheme}, Sci. Record 3(1959), 604--613.

\bibitem{Hager} M. Hager, \emph{Pendant tree-connectivity}, J. Combin.
Theory 38(1985), 179--189.

\bibitem{HAMADA} T. Hamada, T. Nonaka, and I. Yoshimura, \emph{On the Connectivity of Total Graphs}, Math. Ann. 196.30--38 (1972).

\bibitem{Chartrand1}
G. Chartrand, S.F. Kappor, L. Lesniak, D.R. Lick,
\emph{Generalized connectivity in graphs}, Bull. Bombay Math.
Colloq, 2(1984), 1--6.

\bibitem{Chartrand2} G. Chartrand, F. Okamoto, P. Zhang,
\emph{Rainbow trees in graphs and generalized connectivity},
Networks 55(4)(2010), 360--367.

\bibitem{Steeart} G. Chartrand, M. Steeart, \emph{The connectivity of
line graphs}, Math. Ann. 182(1969), 170--174.

\bibitem{Goldsmith3} D.L. Goldsmith, B. Manval, V. Faber,
\emph{Seperation of graphs into three components by removal of
edges}, J. Graph Theory 4(1980), 213--218.

\bibitem{Grotschel1} M. Gr\"{o}tschel,
\emph{The Steiner tree packing problem in $VLSI$ design}, Math.
Program. 78(1997), 265--281.

\bibitem{GLS} R. Gu, X. Li, Y. Shi, The generalized 3-connectivity of random graphs, {\it Acta Math. Sin. (Chin. Ser.)} {\bf 57}(2)(2014), 321--330.


\bibitem{Jain} K. Jain, M. Mahdian, M. Salavatipour,
\emph{Packing Steiner trees}, in: Pro. of the 14th $ACM$-$SIAM$
symposium on Discterte Algorithms, Baltimore, 2003, 266--274.

\bibitem{Kriesell1} M. Kriesell, \emph{Edge-disjoint trees containing
some given vertices in a graph}, J. Combin. Theory, Ser.B 88(2003),
53--65.

\bibitem{Kriesell2} M. Kriesell, \emph{Edge-disjoint Steiner trees
in graphs without large bridges}, J. Graph Theory 62(2009), 188--198.


\bibitem{LM} X. Li, Y. Mao, \emph{Generalized Connectivity of Graphs}, Springer, 2016.

\bibitem{LLSun} H. Li, X. Li, Y. Sun, \emph{The generalied $3$-connectivity
of Cartesian product graphs}, Discrete Math. Theor. Comput. Sci.
14(1)(2012), 43--54.

\bibitem{LLShi} S. Li, X. Li,  Y. Shi, \emph{The minimal size of a graph with generalized connectivity
$\kappa_{3}\geq 2$},  Australasian J. Combin. 51(2011), 209--220.

\bibitem{LST} S. Li, Y. Shi, J. Tu, \emph{The generalized 3-connectivity of Cayley graphs generated
by trees and cycles}, submitted.

\bibitem{LLSS} S. Li, W. Li,  Y. Shi, H. Sun, \emph{On minimally $2$-connected graphs with generalized connectivity $\kappa_{3}=2$}, submitted.

\bibitem{LLL1} S. Li, W. Li, X. Li, \emph{The generalized connectivity of complete
bipartite graphs}, Ars Combin. 104(2012), 65--79.

\bibitem{LLL2} S. Li, W. Li, X. Li, \emph{The generalized connectivity of complete
equipartition $3$-partite graphs}, Bull. Malays. Math. Sci. Soc., in
press.

\bibitem{LLZ} S. Li, X. Li, W. Zhou, \emph{Sharp bounds for the
generalized connectivity $\kappa_3(G)$}, Discrete Math. 310(2010),
2147--2163.

\bibitem{LL} S. Li, X. Li, \emph{Note on the hardness of generalized
connectivity}, J. Combin. Optimization 24(2012), 389--396.

\bibitem{LMS} X. Li, Y. Mao, Y. Sun, \emph{On the generalized (edge-)connectivity
of graphs}, arXiv:1112.0127 [math.CO] 2011.

\bibitem{LMW} X. Li, Y. Mao, L. Wang, \emph{Graphs with large
generalized $3$-edge-connectivity}, arXiv: 1201.3699 [math.CO] 2012.

\bibitem{LYMY} Y. Li, \emph{The generalized 3-(edge) connectivity of line graphs}, submitted.

\bibitem{Nash} C.St.J.A. Nash-Williams,
\emph{Edge-disjonint spanning trees of finite graphs}, J. London
Math. Soc. 36(1961), 445--450.

\bibitem{Oellermann1} O.R. Oellermann, \emph{Connectivity and
edge-connectivity in graphs: A survey}, Congessus Numerantium 116
(1996), 231-252.

\bibitem{Okamoto} F. Okamoto, P. Zhang, \emph{The tree connectivity of
regular complete bipartite graphs}, J. Combin. Math. Combin. Comput.
74(2010), 279--293.

\bibitem{West} D. West, \emph{Introduction to Graph Theory (Second edition)},
Prentice Hall, 2001.

\bibitem{Hamada} T. Hamada and T. Nonaka et al., \emph{On the connectivity of total graph}, Math. Ann. 196.30--38 (1972).

\bibitem{West2} D. West, H. Wu, \emph{Packing Steiner trees and
$S$-connectors in graphs}, J. Combin. Theory, Ser.B 102(2012),
186-205.
\end{thebibliography}
\end{document}